\documentclass{hkart}

\usepackage{amsmath,amsfonts,amssymb,color,enumerate}

\def\nin{\noindent}
\def\defi{:=}
\def\hat{\widehat}
\def\hats #1{\hat{\hat{\hbox{$#1$}}}}
\def\<{\langle} \def\>{\rangle}
\def\.{{\cdot}}

\def\V{{\cal V}}
\def\W{{\cal W}}
\def\AA{{\cal W}\hskip-2pt {\cal A}}
\def\LL{{\cal W}\hskip-1pt {\cal L}}
\def\li #1{#1_{\rm Lie}}
\def\WA{\AA}
\def\A{{\hat G}}

\def\C{\mathbb C}

\def\N{\mathbb N}
\def\P{\mathbb P} 
\def\R{\mathbb R}
\def\RR{\mathcal R}

\def\K{\mathbb K} 
\def\L{\mathfrak L}
\def\T{\mathbb T}
\def\Z{\mathbb Z}

\def\SS{{\mathbb S}}
\def\G{\mathbb G} 
\renewcommand{\Gamma}{{\mathbb G}}
\renewcommand{\Pi}{{\mathbb P}}
\def\GG{\mathcal G}

\def\Span{\mathop{\rm span}\nolimits}

\def\Hom{\mathop{\rm Hom}\nolimits}
\def\Gl{\mathop{\rm Gl}\nolimits}
\def\rank{\mathop{\rm rank}\nolimits}
\def\id{\mathop{\rm id}\nolimits}
\def\pr{\mathop{\rm pr}\nolimits}

\def\inc{\mathop{\rm inc}\nolimits}
\def\UU{\mathop{\bf U\hphantom{}}\nolimits}
\def\bsk{\bigskip}

\def\lead{\leaders\hbox to 1.5ex{\hss${.}$\hss}\hfill}
\def\arr{\hbox to 40pt{\rightarrowfill}}
\def\larr{\hbox to 40pt{\leftarrowfill}}
\def\mapdown#1{\Big\downarrow\rlap{$\vcenter{\hbox{$\scriptstyle#1$}}$}}
\def\lmapdown#1{\llap{$\vcenter{\hbox{$\scriptstyle#1$}}$}\Big\downarrow}
\def\mapright#1{\smash{\mathop{\arr}\limits^{#1}}}
\def\lmapright#1{\smash{\mathop{\arr}\limits_{#1}}}

\def\mapup#1{\Big\uparrow\rlap{$\vcenter{\hbox{$\scriptstyle#1$}}$}}
\def\lmapup#1{\llap{$\vcenter{\hbox{$\scriptstyle#1$}}$}\Big\uparrow}

\newcommand{\End}{\mathrm{End}}
\newcommand{\bL}{\mathbb L}
\newcommand{\bH}{\mathbb H}
\newcommand{\eps}{\varepsilon}
\newcommand{\eK}{{\varepsilon,\K}}

\title{On Weakly Complete Group Algebras\\ of Compact Groups}
                  
\author{Karl Heinrich Hofmann and Linus Kramer\thanks{
Both authors were supported by Mathematisches Forschungsinstitut Oberwolfach
in the program RiP (Research in Pairs). Linus Kramer is
funded by the Deutsche Forschungsgemeinschaft
under Germany's
Excellence Strategy EXC 2044-390685587,
Mathematics M\"unster: Dynamics-Geometry-Structure.}}

\lastname{Hofmann and Kramer}



\address{%
Karl Heinrich Hofmann\\       
Fachbereich Mathematik\\      
Technische Universit\"at Darmstadt\\
Schlossgartenstra{\ss}e 7\\
64289 Darmstadt, Germany\\       
\hbox{hofmann@mathematik.tu-darmstadt.de}}

\address{%
Linus Kramer\\	
Mathematisches Institut\\ 
Universit\"at  M\"unster\\
Einsteinstra{\ss}e 62\\
48149 M\"unster, Germany\\
\hbox{linus.kramer@uni-muenster.de}}

\begin{document}

\maketitle

\vglue-30pt




\begin{abstract} A topological vector space over the real 
or complex field $\K$ is {\it weakly complete} 
if it is isomorphic to a power $\K^J$.
For each topological group $G$ there is a {\it weakly complete 
topological group Hopf algebra} $\K[G]$ over $\K=\R$ or $\C$,
for which three  insights are contributed:
Firstly, {\it there is a comprehensive structure theorem 
saying that the topological algebra $\K[G]$ is
the cartesian product of its finite dimensional minimal ideals 
whose  structure is clarified.}
Secondly, {\it for a  compact {\rm abelian}
group $G$ and its  character group $\A$, the weakly complete 
{\rm complex}
Hopf algebra $\C[G]$ is the product algebra $\C^\A$ with the comultiplication
$c\colon\C^\A\to\C^{\A\times\A}\cong\C^\A\otimes\C^\A$, $c(F)(\chi_1,\chi_2)=
F(\chi_1+\chi_2)$ for $F\colon\A\to\C$ in $\C^\A$. The 
subgroup $\Gamma(\C^\A)$ of grouplike elements of the group of units of the
algebra $\C^\A$ is $\Hom(\A,(\C\setminus\{0\},.))$ 
while
the vector subspace of primitive elements is $\Hom(\A,(\C,+))$.}
This forces  the
group $\Gamma(\R[G])\subseteq\Gamma(\C[G])$ to be 
$\Hom(\A,\SS^1)\cong\hats G\cong G$ with
the complex circle group $\SS^1$. While the relation $\Gamma(\R[G])\cong G$
remains true for {\it any} compact group,  $\Gamma(\C[G])\cong G$ 
holds for a  compact abelian  group $G$ if and only if it is profinite.
Thirdly,  for each pro-Lie algebra $L$ a weakly complete
universal enveloping Hopf algebra $\UU_\K(L)$ over $\K$ exists such
that {\it for each {\em connected} compact group $G$ 
the weakly complete real group Hopf
algebra $\R[G]$ is a quotient Hopf algebra of 
 $\UU_\R(\L(G))$ with the (pro-)Lie algebra
$\L(G)$ of $G$. The group $\Gamma(\UU_\R(\L(G)))$ of grouplike elements
 of the weakly complete enveloping algebra of $\L(G)$ maps onto
$\Gamma(\R[G])\cong G$ and is therefore nontrivial} in   contrast 
to the case of  
the discrete classical enveloping Hopf algebra 
of an abstract Lie algebra.
\end{abstract}

\section{Introduction}
For  much of the material surrounding a theory of group Hopf algebras in the
category of weakly complete real or complex vector spaces we refer
to \cite{dhtwo} and
the forthcoming fourth edition of \cite{compbook}.
Some additional facts are presented here.
A topological vector space over a locally compact field $\K$ is called
{\it weakly complete} if it is isomorphic to $\K^J$ for some set $J$. 
This text considers
$\K=\R$ or $\K=\C$ only and deals with the 
categories $\GG$ of topological groups,
respectively, $\WA$ of weakly complete topological algebras. 
Inside each $\WA$-object $A$
we have the subset $A^{-1}$ of all units (i.e.,  invertible elements) 
which turns out to
be a $\GG$-object in the subspace topology.
In \cite{dhtwo} it was shown that the functor $A\to A^{-1}:\WA\to\G$
has a left adjoint functor 
$G\mapsto \K[G]:\GG\to \WA$. Automatically, $\K[G]$ is 
a topological Hopf algebra. Then for 
each topological group $G$  there is a $\GG$-morphism
$\eta_G\colon G\to\K[G]^{-1}$  such that for each $\WA$-object $A$ and 
$\GG$-morphism $f\colon G\to A^{-1}$ 
there is a unique $\WA$-morphism $f'\colon\K[G]\to A$
such that $f(g)=f'(\eta_G(g))$. The weakly complete algebra
$A=\K[G]$ is seen to have a comultiplication $c\colon A\to A\otimes A$
making it into a symmetric Hopf algebra. The subset 
$$\G(A)=\{a\in A^{-1}:c(a)=a\otimes a\}$$ is a subgroup of $A^{-1}$
and its elements are called {\it grouplike}. 
If $G$ is a discrete group, then $\K[G]$ is the traditional
{\it group algebra} over $\K$ and $\eta_G$ is 
an embedding and its image is $\G(\K[G])$.
In \cite{dhtwo} it was shown that for a compact group $G$ the map
$\eta_G$ is an embedding and for $\K=\R$ induces an isomorphism 
onto $\G(\R[G])$.
We shall see in this text that {\it 
this fails over the complex ground field
 even for all
compact abelian groups which are not profinite.} 
The assignment $G\mapsto\G(\C[G])$ may be viewed as
as `complexification' of the compact group $G$.
This viewpoint becomes important if one considers the module
category of $\C[G]$ in $\W$, i.e. the representations of
$G$ on weakly complete complex vector spaces,  but we will not pursue
this in the present note.

Since we shall deal with compact groups throughout this paper,
we shall always consider $G$ as a subgroup of the group 
$\K[G]^{-1}$ of units of $\K[G]$. In  particular, this means
$G\subseteq \R[G]\subseteq \C[G]$.
 
The article \cite{dhtwo} and the 4th Edition of the book \cite{compbook}
identify a category $\cal H$ 
of weakly complete topological   Hopf algebras for which the functor
$G\to\K[G]$ implements an equivalence of the category of compact groups
and the category~$\cal H$. The vector space continuous dual of $\K[G]$ 
turns out to be the traditional  representation algebra $R(G,\R)$. 
This approach 
 yields a new access to the
Tannaka--Hochschild  duality of the categories of compact groups and
 ``reduced'' real Hopf algebras. 

\nin In all of this, the precise nature of the weakly complete Hopf 
algebras $\K[G]$ even
for compact groups remained somewhat obscure in the nondiscrete case. 
The present paper will present a precise piece of information
on {\it the topological algebra structure of} $\K[G]$ in terms of
a direct product of its finite dimensional minimal ideals whose
precise structure links this presentation with the classical 
information of finite dimensional $G$-modules 
(cf.\ e.g. \cite{compbook}, Chapters 3 and 4). Certain complications
have to be {overcome} on that level 
if one insists {on} an explicit identification
of the algebra and ideal structure of $\R[G]$ as the case $\C[G]$
is easier.
 
 For {\it abelian} compact groups $G$ we shall present
a very direct access to {\it the structure of the weakly complete Hopf
algebra  of} $\C[G]$  by identifying
an isomorphism between $\C[G]$ 
and $\C^{\hat G}$ and by further identifying the
group of grouplike elements to be isomorphic to $(\L(G),+)\oplus G$ with
the pro-Lie algebra $\L(G)$ of $G$ (see \cite{compbook}).
 The subgroup $G$ of $\C[G]^{-1}$ therefore agrees with the
group of grouplike elements of $\C[G]$ 
if and only if $\L(G)=\{0\}$ if and only if $G$ is totally disconnected.

\medskip

\nin The presence of weakly complete Lie algebras over $\K$ 
in the group $\K$-Hopf algebra of a 
compact group motivates a proof of the  existence of 
a weakly complete  universal enveloping
algebra over $\K$ for $\WA$ Lie algebras over $\K$. 
In  contrast to the case of enveloping algebras
on the purely algebraic side, the $\WA$ enveloping Hopf algebras 
will  sometimes  have grouplike
elements. The universal property of the $\WA$ enveloping 
Hopf  algebras will show that
each $\WA$ Hopf-group algebra $\K[G]$ of a 
compact {\it connected} group $G$ is a quotient algebra
of the  $\WA$ enveloping algebra of $\L(G)$. So
$\WA$ enveloping algebras have a tendency 
of being  larger than 
$\WA$ group algebras. 

For a special class of profinite dimensional
Lie algebras a similar but different approach to appropriate enveloping 
algebras is considered in \cite{ham}.

A preprint of the present material
appeared in the Series of Preprints of the 
Mathematical Research Institute of
Oberwolfach \cite{hofkra}. 
 
\medskip

\section{Weakly Complete Hopf Algebras}

For the basic theory of weakly complete Hopf algebras we may safely
refer to \cite{dhtwo} and \cite{compbook}, 4th Edition. For the present 
discussion we need a reminder of some basic concepts. 

\begin{Definition} Let $A$ be a weakly complete symmetric Hopf algebra,
i.e.\ a group object in the monoidal category $(\W,\otimes_W)$ of
weakly complete vector spaces (see \cite{compbook}, Appendix 7 and
Definition A3.62), with
comultiplication $c\colon A\to A\otimes A$ and coidentity $k\colon A\to \K$. 

 An element $a\in A$ is called
{\it grouplike} if $c(a)=a\otimes a$ and $k(a)=1$. The subgroup
of grouplike elements in the group of units $A^{-1}$ will be denoted
$\G(A)$. 

An element  $a\in A$ is called {\it primitive}, if
$c(a)= a\otimes 1 + 1\otimes a$. The Lie algebra of primitive elements of 
$\li A$, i.e.\ the weakly complete Lie algebra obtained by endowing the
weakly complete vector space underlying $A$ with the Lie bracket
obtained by $[a,b]=ab-ba$, will be
denoted $\Pi(A)$.                                                              
\end{Definition}

\medskip

Any weakly complete
symmetric Hopf algebra $A$  has an
exponential function $\exp_A\colon A\to A^{-1}$
as explained in \cite{dhtwo}, Theorem 3.12 or in \cite{compbook},
4th Edition, A7.41.

\begin{Theorem} \label{exp}
 Let $A$ be a weakly complete symmetric Hopf algebra.
Then the following statements hold:
\begin{enumerate}[\rm(i)]
\item The  set $\Gamma(A)$
of  grouplike elements of a weakly complete
symmetric Hopf algebra $A$ is a closed subgroup of $(A,{\cdot})$
and therefore is a pro-Lie group.
\item The set $\Pi(A)$ of  primitive elements
of $A$ is a closed Lie subalgebra of $A_{\rm Lie}$ 
 and therefore is a pro-Lie algebra.
\item $\Pi(A)\cong\L(\G(A))$ and 
the exponential function $\exp_A$ of $A$  induces the exponential
function $\exp_{\Gamma(A)}\colon \Pi(A)\to \Gamma(A)$
of the pro-Lie group $\Gamma(A)$.
\end{enumerate}
\end{Theorem}
For a proof  see  e.g.\ \cite{dhtwo}, Theorem 6.15.

\begin{Definition} 
For an arbitrary topological group $G$ we define
$R(G,\K)\subseteq C(G,\K)$ to be that set of continuous functions
$f\colon G\to \K$ for which the linear span
of the set of translations ${}_gf$, ${}_gf(h)=f(hg)$,
 is a finite dimensional vector subspace of $C(G,\K)$.
The functions in $R(G,\K)$ are called {\it representative
functions.} \end{Definition}

Clearly $R(G,\K)$ is a subalgebra of $C(G,\K)$ also
known as the {\it representation algebra} of $G$.
In \cite{dhtwo}, Theorem 7.7(a) the following {duality} result was
{shown.}

\begin{Theorem} \label{4.7} 
{\rm(The Dual of a Weakly Complete Group Algebra  $\K[G]$)}\\
  For an arbitrary topological group $G$,
the function
$$ F_G\colon \K[G]'\to \RR(G,\K),\quad F_G(\omega)=\omega\circ \eta_G$$
is a natural isomorphism of Hopf algebras.
\end{Theorem}

This applies, of course, to compact groups, in which case the
Hopf algebra $R(G,\K)$ is a well-known object.

For easy reference we record the following facts in the case of
a compact group $G$ for which we recall
$G\subseteq \K[G]$:

\begin{Theorem} \label{3.2 real} 
For any {\rm compact} topological group $G$, the following statements hold:
\begin{enumerate}[\rm(i)]
  \item We have
$G\subseteq \Gamma(\K(G))\subseteq\K[G]^{-1}$

\item In the case of  $\K=\R$ the
equality $G=\Gamma(\R[G])$ holds.
\end{enumerate}
\end{Theorem}
For (1) see \cite{dhtwo}, 5.4, and for (ii) see
\cite{dhtwo}, 8.7.

\noindent
For the complex case, we shall see later in this paper
that in many cases, a compact group $G$ is a
 {\em proper} subgroup of $\Gamma(\C[G])$.

\section{Some Preservation Properties of $\K[-]$}

Let us explicitly formulate and prove some preservation properties of 
our functor $\K[-]$. Left adjoint functors preserve epics. A morphism
of compact groups is an epimorphism if and only if it is surjective
(see \cite{compbook}, RA3.17). Therefore the following lemma is to be
expected. 

\begin{Lemma}  For every surjective morphism 
$f\colon G\to H$ of {\it compact} groups
the morphism $\K[f]\colon \K[G]\to \K[H]$ of weakly complete 
$\K$- Hopf algebras is  surjective.
\end{Lemma}

\begin{Proof} From the surjectivity of $f\colon G\to H$ we conclude
that 

\centerline{$f(\Span(G)) = \Span(f(G)) =\Span(H)$}

\nin is dense in $\K[H]$ by
Proposition 5.3 of \cite{dhtwo}, and likewise $\K(f)(\K[G])$
is dense in $\K[H]$. But $\K[f]$ is, in particular, a $\W$-morphism,
that is, a morphism of weakly complete vector spaces. 
Every such has a closed image by  
\cite{compbook}, Theorem 7.30(iv). Hence $\K(f)(\K[G])=\K[H]$.
\end{Proof}
 
This particular left adjoint  functor $\K[-]$, however, also preserves
the injectivity of morphisms:
\medskip

\begin{Theorem} \label{7.2}  If $G$ is a closed subgroup of 
the {\rm compact} group $H$,
then $\K[G]\subseteq \K[H]$ (up to natural isomorphism).
\end{Theorem}

\begin{Proof} From the injectivity of a morphism of compact groups
$j\colon G\to H$ we derive the surjectivity of
$C(j,\K)\colon  C(H,\K)\to C(G,\K)$ by the Tietze Extension Theorem.
Now we set $M\defi C(f,\K)\big(R(H,\K))\subseteq R(G,\K)$. 
Since $\R(H,\K)$ is dense in $C(H,\K)$ in the norm topology,
$M$ is dense in $R(G,\K)$ in the norm topology. Then it is dense
in $L^2(G,\K)$ in the $L^2$-topology, and $M$ is a $G$-module.
In the case of $\K=\R$ we can now apply Lemma 8.11 of \cite{dhtwo} and
conclude that $M=R(G,\R)$. Thus $R(G,j)\colon R(H,\R)\to R(G,\R)$
is surjective. By Theorem 7.7 of \cite{dhtwo} this implies that
$\R[f]'\colon \R[H]'\to\R[G]'$ is surjective. The duality between
$\K$-vector spaces and weakly complete $\K$-vector spaces shows that
$\R[f]\colon\R[G]\to\R[H]$ is injective. This proves the theorem for
$\K=\R$.
But then the commuting diagram
\medskip

$$\begin{matrix} \C\otimes\R[G]&\mapright{\C\otimes\R[f]}&\C\otimes\R[H]\\                                   
          \lmapdown{\cong}&&\mapdown{\cong}\\                                                         
           \C[G]&\lmapright{\C[f]}&\C[H]\end{matrix}
$$
\medskip

\nin shows that $\C[f]$ is also injective.

In the category of weakly complete vector spaces every injective
morphism is an embedding by duality   since every surjective morphism of vector
spaces is a coretraction. \end{Proof}

\begin{Corollary} \label{7.3} Let $G_0$ denote the identity component of
the compact group $G$. Then
\begin{enumerate}[\rm(i)]  
\item The Hopf algebra $\K[G_0]$ is a Hopf subalgebra of $\K[G]$.
\item  $\R[G_0]$ is algebraically and topologically generated by 
$\Pi(\R[G])\cong\L(G)$.
\end{enumerate}
\end{Corollary}

\begin{Proof} (i) is a consequence of Theorem \ref{7.2}.

 (ii)  The compact group $G_0$ is algebraically and topologically
generated by $\exp_G(\L(G))$ (cf.\ \cite{probook}, Corollary 4.22, p.~191),
and $\overline{\Span(G_0)}=\R[G_0]$ by \cite{dhtwo},
Corollary 5.3.
\end{Proof} 

\section{A Principal Structure Theorem of $\K[G]$}

Let $G$ be a compact group and let $E$ be a finite dimensional vector 
space over $\K\in\{\R,\C\}$. 
We recall that the {\it character of a representation}
$\rho\colon G\to \mathrm{End}_\K(E)$ is the continuous 
map $g\mapsto\mathrm{tr}_\K(\rho(g))$. We also say that $\chi$ is the
character of the $G$-module $E$.
A representation, respectively, a $G$ module, 
is determined by its character  up to isomorphism.
A character is called {\it irreducible} if the corresponding representation is
irreducible (over the ground field $\K)$, equivalently, the corresponding
$G$-module is simple. We denote the 
set of all irreducible characters of $G$ over $\K$ by $\hat G_\K$.
For each character $\chi$ of $G$, we select a 
finite dimensional $G$-module $E_{\chi,\K}$ having $\chi$ as its character.
If $\eps$ is an irreducible character, then the ring $\bL_\eK=\End_G(E_\eK)$
of all $\K$-linear endomorphisms of $E_\eK$ which commute with the 
$G$-action
is, by Schur's Lemma, a finite dimensional division ring over $\K$.
Hence $$\begin{matrix}
\bL_\eK&=&\C\hfill&\mbox{if}&\K=\C,\\
\bL_\eK&\in&\{\R,\C,\bH\}&\mbox{if}&\K=\R,\\
\end{matrix}$$
where $\bH$, as is usual, denotes the skew-field of quaternions.
We view $E_\eK$ as a right module over $\bL_\eK$.
We denote the corresponding representation by
$
\rho_\eK:G\to\End_{\bL_\eK}(E_\eK)\subseteq\End_\K(E_\eK).
$
Before we enter the presentation of the principal theorem on the
weakly complete group algebra $\K[G]$ of a compact group we 
 elaborate on some basic ideas of finite dimensional representation
theory, indeed extending some of the presentation such as 
it can be found e.g.\
in Chapter 3 of \cite{compbook}. The first lemma extends the details of
Proposition 3.21 of \cite{compbook} and the comments which precede it.
\begin{Lemma} \label{DensityCor}
Let $E$ be a finite dimensional vector space over $\K$ and 
$\rho\colon G\longrightarrow\End_\K(E)$
an irreducible representation of a group $G$.
Let $A$ denote the $\K$-span of the set $\{\rho(g)\mid g\in G\}$.
Then $A=\End_\bL(E)$, where $\bL=\End_A(E)=\End_G(E)$, 
$\bL\in\{\R,\C,\bH\}$.
\end{Lemma}
\begin{proof}
First of all we note that $A$ is a $\K$-algebra containing $\id_E$.
Hence every $A$-submodule of the additive group $E$ is a 
$G$-invariant linear subspace, and vice versa.
Therefore $E$ is a simple $A$-module and so 
Jacobson's Density Theorem  applies, which, for the sake of
completeness, we cite here in its entirety (see e.g.\ \cite{CurtisReiner}).

\smallskip

\noindent {\bf Jacobson's Density Theorem.} \label{density}
{\it 
Let $M\neq 0$ be an (additive) abelian group, let 
$A\subseteq\End(M)$ be a subring and suppose that $M$ is
simple as a left $A$-module. Put $\bL=\End_A(M)$. Then $\bL$ is a division ring
and $M$ is a right $\bL$-module in a natural  way.
For every $2k$-tuple $(x_1,\ldots,x_k,y_1,\ldots,y_k)\in M^{2k}$, such that the 
elements $x_1,\ldots,x_k$ are linearly independent, 
there exists $a\in A$ such that 
$a(x_i)=y_i$ holds for all $i=1,\ldots k$.}
\smallskip

\noindent
Now the division ring $\bL$ is a finite dimensional $\K$-algebra over $\K$, 
and hence is
isomorphic to $\R,\C$, or $\bH$. Moreover, $A\subseteq\End_\bL(E)$.
Let $x_1,\ldots,x_m$ be a $\bL$-basis for 
$E$, and let $\phi\in\End_\bL(E)$ be arbitrary. 
Then there exists an element $a\in A$
such that $a(x_i)=\phi(x_i)$ holds for all $i=1,\ldots,m$.
Therefore $\phi=a$ and thus $\End_\bL(E)=A$.
\end{proof}

\nin The following result now extends \cite{compbook},~Lemma 3.14.
\begin{Lemma} \label{HomCor}
Let $E$ and $F$ be finite dimensional vector spaces over $\K$ 
and suppose that 
$\rho\colon G\to\End_\K(E)$ and $\sigma\colon H\to\End_\K(F)$ 
are irreducible representations of groups $G,H$.
Suppose also that $\End_G(E)=\bL=\End_H(F)$. 
Then $\Hom_\bL(F,E)$ is an irreducible
$G\times H$-module over $\K$, 
where $(g,h)(f)=\rho(g)\circ f\circ\sigma(h^{-1})$.
\end{Lemma}
\begin{proof}
We define $A\subseteq\End_\K(E)$ and $B\subseteq\End_\K(F)$ 
as in Corollary~\ref{DensityCor}.
Then $A=\End_\bL(E)$  and $B=\End_\bL(F)$.
The $\K$-vector space $\Hom(F,E)$ is in a natural way a right 
$A$-module and a left $B$-module.
For every nonzero $f\in\Hom_\bL(F,E)$ we have $AfB=\Hom_\bL(F,E)$.
Therefore $\Hom_\bL(F,E)$ is simple as $G\times H$-module over $\K$.
\end{proof}

We are now ready to prove a principal structure theorem 
{for} the weakly complete
group algebra $\K[G]$ of a compact group $G$ for either $\K=\R$ or
$\K=\C$. 
For each $\eps\in\hat G_\K$ we have the $G$-module $E_{\eps,\K}$ and
the corresponding irreducible representation 
$\rho_{\eps,\K}\colon G\to\End_{\eps,\K}(E_{\eps,\K})$ into the group of
units of the concrete {\it matrix ring}
$M_\eps\defi\End_{\eps,\K}(E_{\eps,\K}))$ over $\bL=\bL_{\eps,\K}$
of $\bL$-dimension {$(\dim_\bL E_{\eps,\K})^2$}. 
Accordingly there is a unique function
$\rho_G\colon G\to \prod_{\eps\in\hat G_\K} M_\eps$ which is an injective 
group morphism into the multiplicative group of units of the product
defined by the universal property of the product such that
$$\begin{matrix} 
G&\mapright{\rho_G}&\prod_{\eps\in\hat G_\K}M_\eps\\
\lmapdown{=}&&\mapdown{\pr_\chi}\\
G&\lmapright{\rho_{\chi,\K}}&M_\chi\\
\end{matrix}$$
commutes for all $\chi\in\hat G_\K$.

\begin{Theorem} \label{principal thm}
For any compact group $G$ the
weakly complete symmetric Hopf algebra $\K[G]$
is a direct product
$$ \K[G]=\prod_{\epsilon \in\hat G_\K}\K_\epsilon[G]$$
of finite dimensional minimal two-sided ideals 
$\K_\epsilon[G]$ such that  for each $\eps\in\hat G_\K$ 
there is a $\K$-algebra isomorphism
$$
\K_\eps[G]\cong M_\eps=\End_{\bL_\eK}(E_{\eps,\K}).
$$
In particular, each of these two-sided ideals $\K_\eps[G]$ is a 
two sided simple $G\times G$-module  
and as an algebra is isomorphic to a full matrix ring over $\bL$.
\end{Theorem} 
\begin{Remark}
The diagram
$$\begin{matrix}
G&\mapright{\eta_G}& \K[G] \\
\lmapdown{=} &&\mapdown{\cong}\\
G& \lmapright{\rho_G} & \prod_{\eps\in\hat G_\K}M_\eps\\
\lmapdown{=}&&\mapdown{\pr_\chi}\\
G&\lmapright{\rho_{\chi,\K}}&M_\chi\\
\end{matrix}$$
commutes for all $\chi\in\hat G_\K$.
\end{Remark}

\begin{proof}
By Theorem 2.4 and \cite{compbook} Theorem 3.28,  the topological
dual $\K[G]'\cong R(G,\K)$ is the direct sum of the finite dimensional
two-sided $G$-submodules $R_\epsilon(G,\K)$ as $\epsilon$ ranges through
the set of irreducible characters in $\hat G_\K$.
The $G\times G$-module $R_\epsilon(G,\K)$ is defined in \cite{compbook}
as the image of the linear map 
$$
\phi:E_\eK'\otimes_\K E_\eK\longrightarrow R(G,\K),
$$ 
where $$\phi(u\otimes v)(g)=\langle u,\rho_\eK(g)v\rangle.$$
If we put $\psi(f)(g)=\mathrm{tr}_\K(\rho_\eK(g)f)$, for 
$f\in\End_\K(E_\eK)$ and $g\in G$, then the diagram
$$
\begin{matrix}
E'_\eK\otimes_\K E_\eK &\mapright{\phi} & R_\eps(G,\K) \\
\lmapdown{s} &&\lmapdown{=}\\
\End_\K(E_\eK)&\mapright{\psi} & R_\eps(G,\K)
\end{matrix}
$$
commutes, where $s(u\otimes v)=[w\mapsto v\langle u,w\rangle]$.
We recall that group $G\times G$ acts on 
$R_\eps(G,\K)$ via $(a,b)(\lambda)=[g\mapsto \lambda(a^{-1}gb)]$.
If we put $(a,b)(u\otimes v)=(u\circ\rho_\eK(a^{-1}))\otimes\rho_\eK(b)v$
and $(a,b)(f)=\rho_\eK(b)\circ f\circ\rho_\eK(a^{-1})$, then
all maps in this diagram are $G\times G$-equivariant.

Suppose that $\K=\bL$. Then $\End_\K(E_\eK)=\End_L(E_\eK)$ is simple
as a $G\times G$-module by Lemma~\ref{HomCor} 
and thus $\psi$ is an isomorphism.

Suppose next that $\K\subsetneq\bL$. Then $\K=\R$ and 
$\bL=\C$ or $\bL=\bH$. 
By the averaging process in \cite{compbook} Lemma~2.15
there exists a $G$-invariant positive definite $\bL$-hermitian form
$(\cdot |\cdot)$ on $E$, semilinear in the first argument and linear in the second argument.
This allows us to rewrite $R_\eps(G,\K)$
as the span of the maps $g\mapsto \mathrm{Re}(u|gv)$, for $u,v\in E$.
The $G$-invariance of $(\cdot|\cdot)$ yields that 
$\mathrm{Re}(au|gbv)=\mathrm{Re}(u|a^{-1}gbv)$.
If we consider the algebra inclusion 
$$j\colon\End_\bL(E_\eK)\to \End_\K(E_\eK),$$
then $\mathrm{Re}(u|v)=\mathrm{tr}_\K[w\mapsto v(u|w)]$ 
holds for the trace map of $\End_\K(E_\eK)$.
It follows that the map $\psi\circ j$ in the diagram
$$
\begin{matrix}
E'_\eK\otimes_\K E_\eK &\mapright{\phi} & R_\eps(G,\K) \\
\lmapdown{s} &&\lmapdown{=}\\
\End_\eK(E_\eK)&\mapright{\psi} & R_\eps(G,\K)\\
\lmapup{j} && \\
\End_\bL(E_\eK)
\end{matrix}
$$
is surjective and $G\times G$-equivariant.
Since $\End_\bL(E_\eK)$ is a simple $G\times G$-module over $\K$
by Lemma~\ref{HomCor}, the map $\psi\circ j$ is an isomorphism.

\hskip -15pt For the remaining part of the proof we apply standard duality
theory.
We put $$R^\chi=\bigoplus_{\chi\neq\eps\in\hat G_\K}R_\eps(G,\K)$$
and define $\K_\chi[G]$ as the annihilator of $R^\chi$.
The annihilator mechanism supplies us with the diagram
$$\begin{matrix}%
\K[G]&{\buildrel\perp\over\longleftrightarrow}&\{0\}&\cr
   \Big|&&   \Big|&\cr
\K_\chi[G]&{\buildrel\perp\over\longleftrightarrow}&R^\chi&\cr
   \Big|&&   \Big|&\Big\}\cong R_\chi(G,\K)\cr
\{0\} &{\buildrel\perp\over\longleftrightarrow}&R(G,\K).&\cr
\end{matrix}$$

\smallskip\nin
By the duality of $\V_\K$ and $\W_\K$ it follows that
$\K[G]\cong\prod_{\epsilon\in{\hat G}_K}\K_\epsilon[G]$
with $$\K_\epsilon[G]\cong R_\epsilon(G,\K)'.$$

Now, if  any closed vector subspace $J$ of  $\K[G]$
satisfies $G\.J\subseteq J$ and $J\.G\subseteq J$,
then we also have
$\Span(G)\.J\subseteq J$ and $J\.\Span(G)\subseteq J$
(where we view $G$ as a subset of $\K[G]$). 
Then  Proposition 5.3 of \cite{dhtwo} says that $\overline{\Span(G)}=\K[G]$,
 and so $\K[G]\.J\subseteq J$ and $J\.\K[G]\subseteq J$. That is,
$J$ is a closed two-sided ideal of $\K[G]$. Therefore each $\K_\epsilon[G]$
is a two-sided ideal in $\K[G]$.

It remains to clarify the multiplicative structure 
of the ideals $\K_\eps[G]$.
If we consider $\eps\in{\hat G}_\K$ 
and the representation $\rho_\eK$,
then the map 
$$G\mapright{\rho_\eK}\Gl_{\bL_\eK}(E_{\eps,\K})=\End_{\bL_\eK}(E_{\eps,\K})^{-1}\mapright{\inc} \End_{\bL_\eK}(E_{\eps,\K})$$
and the universal property of $\K[G]$ described in the Weakly Complete Group
Algebra Theorem 5.1 of \cite{dhtwo}
provides a morphism of weakly complete algebras 
$$\pi_\eps\colon\K[G]\to\End_{\bL_\eK}(E_{\eps,\K})$$ extending $\rho_\eK$. 
We also have the product projection of weakly complete
algebras $\pr_\epsilon\colon \K[G]\to \K_\epsilon[G]$. 
Both maps $\pi_\eps$ and $\pr_\epsilon$ have the same
kernel $\prod_{\epsilon\ne\epsilon'\in{\hat G}_\K}A_{\epsilon'}$.
So there is an injective morphism 
$$
\alpha\colon \K_\epsilon[G]\to \End_{\bL_\eK}(E_\eK)
$$
such that $\pi_\eps=\alpha\circ\pr_\epsilon$.
Since both algebras have the same dimension, 
$\alpha$ is an isomorphism of $\K$-algebras.
\end{proof}

\begin{Corollary}
There is an isomorphism of $G\times G$-modules
$$
R(G,\K)=\bigoplus_{\eps\in\hat G_\K}\End_{\bL_\eK}(E_\eK).
$$
Thus the multiplicity $m$ of $E_\eK$ as a $G$-module in $R(G,\K)$
is $$m=\dim_{\bL_\eK}(E_{\eK})=\frac{\dim_\K E_{\eps,\K}}{\dim_\K \bL}.$$
\end{Corollary}
This conclusion  is well-known for $\K=\C$ 
(see e.g Theorems 3.22 and 3.28 in \cite{compbook})
 but we could not readily find a reference for $\K=\R$.
While the algebra structure of the weakly complete 
symmetric Hopf algebra
$\K[G]$ is satisfactorily elucidated in Theorem 
\ref{principal thm}, the comultiplication seems to
be not easily accessible due to complications of
the way how the representation theory of $G\times G$
reduces to that of $G$ in general. In the case 
of commutative compact groups $G$ and the complex 
ground field $\C$ these complications go away, and so
we shall clarify the situation in these circumstances
in the subsequent section.

\section{The Weakly Complete Group Algebras\\
of Compact Abelian Groups: 
An Alternative View}
We have seen the usefulness of the concept of a weakly complete
group algebra $\K[G]$ over the real or complex numbers. We obtained
its existence from the Adjoint Functor Existence Theorem. 
This is rather remote
from a concrete construction. 
It may therefore be helpful to see the whole apparatus
in a much more concrete way at least for a substantial subcategory of
the category of compact groups, namely, 
the category of compact abelian groups
for which we already have a familiar duality theory due to {\sc Pontryagin}
and {\sc Van Kampen} (see e.g.\  \cite{compbook}, Chapter 7).

In this section let $G$ be a compact abelian group and
$\hat G={\mathcal{C\hskip-2pt AB}}(G,\T)$ 
(with the category $\mathcal{C\hskip-2pt AB}$ of compact abelian groups
and $\T=\R/\Z$)
its discrete character group.
These groups are written additively. For $\K=\C$ there is a natural
bijection $\hat G\to\hat G_\C$ from the character group to the set
of equivalence classes of complex simple $G$-modules (cf.\ 
\cite{compbook}, Lemma 2.30 (p.43), Exercise E3.10 (p.66)), and also 
Proposition 3.56 (p.87)
for some information on $\hat G_\R$). 
This bijection associates with a character $\chi\in\hat G=
\Hom(G,\T)$ the class of the module $E_\chi=\C$, 
$\chi\.c=e^{2\pi i\chi}c$.  Accordingly, \cite{compbook}
Theorem 3.28 (12) reads
$R(G,\C)= \sum_{\chi\in\hat G} \C\.f_\chi$, for a suitable basis 
$f_\chi$, $\chi\in\hat G$, $f_\chi(g)=e^{2\pi i\chi(g)}$. 
In other words, as a $G$-module, $R(G,\C)\cong \C^{(\hat G)}$.
Accordingly, we expect $\C[G]$ to be uncomplicated. Our 
Theorem \ref{principal thm} makes this clear: 

\smallskip\noindent
{\it The complex algebra $\C[G]$ may be
naturally identified with the componentwise algebra $\C^{\hat G}$.}

\smallskip

In the abelian case, our understanding of the comultiplication
of $\C[G]=\C^{\hat G}$ is much more explicit than in the general
situation of Theorem \ref{principal thm}.
 Each character $\chi\colon G\to \T$
determines a morphism $f_\chi\colon G\to \C^{-1}=\C^\times$, 
$f_\chi(g)=e^{2\pi i\chi(g)}$, $g\in G\subseteq \C^{\hat G}$. 
By the universal property of 
$\C[G]=\C^{\hat G}$, this value agrees with the $\chi$-th
projection of $g\in G\subseteq \C^{\hat G}$. Hence
$$(\forall g\in G, \chi\in\hat G)\,\eta_G(g)(\chi)
=e^{2\pi i\<\chi,g\>}.$$
Accordingly, if we write $\SS^1=\{z\in \C;|z|=1\}$, then
$g\in \Hom(\hat G,\SS^1)\cong\hats G\cong G$.
Then in view of $G\subseteq \R[G]\subseteq\C[G]$
we  have
$$ \Hom(\hat G,\SS^1)\subseteq \overline{\Span_\R(\Hom(\hat G,\SS^1))}
=\R[G]\subseteq\C[G]=\C^{\hat G}.$$

\medskip

 Recall from \cite{dhtwo}, Theorem 5.5 that we have an isomorphism
$$\alpha_G\colon \C[G\times G]\to \C[G]\otimes_\W \C[G],$$ and from
\cite{dhtwo} Lemma 5.12 we recall the comultiplication 
$\gamma_G\colon\C[G]\to \C[G]\otimes_\W\C[G]$ to be the composition
$$\C[G]\mapright{\delta_G}\C[G\times G]
  \mapright{\alpha_G}\C[G]\otimes_\W \C[G].$$
Now for a compact abelian group $G$,
 the diagonal morphism $\delta_G\colon G\to G\times G$ has the
group operation of $\hat G$ as its dual, namely: 
$$\hat{\delta_G}\colon \hat G \times \hat G\to\hat G.\quad 
\hat{\delta_G} (\chi_1,\chi_2)=\chi_1+\chi_2,$$
as we write abelian group operations additively in general.
If now we  also write
$\C[G]\otimes_\W \C[G] =\C^{\hat G\times \hat G}$ (identifying 
$\phi\otimes \psi$ with $(\chi_1,\chi_2)\mapsto \phi(\chi_1)\psi(\chi_2)$),
then we have
$$\gamma_G=\C^{\hat{\delta_G}}\colon \C^{\hat G}\to\C^{{\hat G}\times{\hat G}},
\hbox{i.e., }
(\forall \phi\in\C^{\hat G}),  
\gamma_G(\phi)(\chi_1,\chi_2)=\phi(\chi_1+\chi_2).
$$
This allows us to determine explicitly the elements of the
group $\G(\C^{\hat G})$ of all grouplike elements: 

Indeed a nonzero
element $\phi\in \C^{\hat G}$ is in $\G(\C^{\hat G})$ if and only if
$$ \gamma_G(\phi)= \phi\otimes \phi\hbox{\quad in\quad  
  $\C^{\hat G}\otimes_\W \C^{\hat G}=\C^{\hat G\times\hat G}$},$$  
    where $(\phi\otimes\phi)(\chi_1,\chi_2)=\phi(\chi_1)\phi(\chi_2)$. 
This is the case if and only if 
$$(\forall \phi_1, \phi_2\in\hat G)\, 
\phi(\chi_1+\chi_2)=\gamma_G(\phi)(\chi_1,\chi_2)
=(\phi\otimes\phi)(\chi_1,\chi_2)=\phi(\chi_1)\phi(\chi_2),$$
that is, if and only if $\phi$ is a morphism of groups from $\hat G$ to 
$\C^\times=(\C\setminus\{0\}, \cdot)$.

\medskip

 Similarly, an element $\phi\in\C^{\hat G}$ is primitive if and only if
$$\phi(\chi_1+\chi_2)=\gamma_G(\phi)(\chi_1,\chi_2)=
\big((\phi\otimes1)+1\otimes\phi)\big)(\chi_1,\chi_2)
=\phi(\chi_1)+\phi(\chi_2)$$
if and only if $\phi\colon \hat G\to (\C,+)$ is a morphism of topological 
groups.

Let us summarize this discourse:

\begin{Theorem} \label{abelian thm} 
{\rm (The Group Hopf Algebra $\C[G]$ for
 Compact Abelian $G$)} \quad
Let $G$ be a compact abelian group and $A$ its weakly complete
commutative symmetric group Hopf algebra $\C[G]$ and let 
$\hat G=\Hom(G,\T)$
be its character.
\begin{enumerate}[\rm(i)]
\item Then $A$ may be
identified with $\C^{\hat G}$ 
such that $g\colon \hat G\to A^{-1}$ is defined by
$$ (\forall  \chi\in\hat G)\,g(\chi)=e^{2\pi i\<\chi,g\>}\in\SS^1,$$
where $\SS^1=\{z\in\C:|z|=1\}\subseteq \C^\times$.
The natural image of $G$ in $A^{-1}$ is
$$G=\Hom(\hat G,\SS^1)\cong\hats G,$$ 
and
$$G=\Hom(\hat G,\SS^1)\subseteq \R[G]\subseteq \C[G]=\C^{\hat G}.$$

\item If, as is possible in the category of weakly complete
vector spaces, the weakly complete vector spaces 
$A\otimes_\W A$ and $\C^{\hat G\times\hat G}$ are identified,
then the comultiplication 
$\gamma_G\colon A\to A\otimes_\W A$ of $A$ is given by
$$(\forall \phi\colon \hat G\to \C,\,  \chi_1,\chi_2\in\hat G)\quad  
\gamma_G(\phi)(\chi_1,\chi_2)= \phi(\chi_1+\chi_2)\in\C.$$

\item The group of grouplike
elements of $A$ is 
$$\G(A)=\Hom(\hat G,\C^\times)\subseteq \C^{\hat G}.$$

\item The weakly complete Lie algebra of primitive elements of
$A$ is 
$$\P(A)=\Hom(\hat G,\C)\subseteq \C^{\hat G}.$$
\end{enumerate}
\vskip-30pt \qed
\end{Theorem}

\noindent
We write $\R^\times_+$ for the multiplicative
subgroup  $\{z\in\C: 0<z\in\R\subseteq\C\}$ of $\C^\times$.

\begin{Corollary} \label{abelian primitive} 
For a compact abelian group $G$ and the weakly complete
commutative unital algebra $A\defi\C[G]$ we have a commutative
diagram
$$\begin{matrix}\P(A)&=&\Hom(\hat G,\R)+\Hom(\hat G,i\R)&
                              \mapright{\cong}&\L(G)\times\L(G)\cr
     \lmapdown{\exp_A} &&&& \mapdown{\id_{\L(G)}\times \exp_G}\cr
          \G(A)&=&\Hom(\hat G,\R^\times_+)\.\Hom(\hat G,\SS^1)&
                     \lmapright{\cong}&\hskip-15pt\L(G)\times G.\cr
\end{matrix}$$
The unique maximal compact subgroup of $\G(A)$ is 
$G=\Hom(\hat G,\SS^1)$. 
\end{Corollary}
\begin{Proof} There is an elementary isomorphism of topological groups
$$(r,t+\Z)\mapsto e^re^{2\pi it}=e^{r+2\pi it}:\R\times \T\to \C^\times.$$
 Accordingly, \vskip-17pt 
$$\G(A)=\Hom(\hat G,\C^\times)\cong \Hom(\hat G,\R)\oplus \Hom(\hat G,\T).$$
Now $\Hom(\hat G,\R)\cong\Hom(\R,G)$ (cf.\ \cite{compbook},
 Proposition 7.11(iii)),
 $\Hom(\R,G){=}\L(G)$ by \cite{compbook}, Definition 5.7 
(cf.\ Proposition 7.36ff., Theorem 7.66)
and $\Hom(\hat G,\T){=}\hats G\cong G$ by \cite{compbook}, Theorem 2.32.
For the exponential function $\exp_A$ of a weakly complete unital 
symmetric Hopf algebra is treated in Theorem 2.2 above.
\end{Proof}

A compact abelian group is totally disconnected (i.e.\ profinite) if and only if
$\L(G)=\{0\}$ (cf. \cite{compbook}, Corollary 7.72).

\begin{Remark} For a compact abelian group $G$ the equality
$G=\G(\C[G])$ holds if and only if $G$ is totally disconnected (i.e.\ profinite).
\end{Remark}
\begin{Proof} By Theorem \ref{abelian thm}, the equality holds 
if and only if $L(G)=\{0\}$ if and only if 
$\Hom(\hat G,\R)=\{0\}$ if and only if $\hat G$ is a torsion group
(cf.\ e.g.\ \cite{compbook}, Propositions A1.33, A1.39) if and only if $G$ is totally
disconnected (see Corollary 8.5).
\end{Proof}
\noindent In particular, e.g., $\T\ne \G(\C[\T])$.

\medskip

Now we understand $\C[G]=\C^{\hat G}$ rather explicitly, but  $\R[G]$ 
only rather implicitly. However, Theorem \ref{principal thm} applies with
$\K=\R$ in order to shed some light on its intrinsic structure.

We define the function $\sigma_G\colon \C[G]\to\C[G]$ 
as follows: For $\chi\in\hat G$ we set $\check \chi(g)=\chi(-g)=-\chi(g)$.
Then we define
$$(\forall \phi\in\C^{\hat G})\, 
\sigma(\phi)(\chi)=\overline \phi(\check \chi).$$

\begin{Exercise} \label{involution} 
For a compact abelian group $G$, the function
$\sigma_G$ is an involution of weakly complete real algebras of $\C[G]$
whose precise fixed point algebra is $\R[G]$. 
Accordingly, $\C[G]=\R[G]\oplus i\R[G]$.
\end{Exercise}

\begin{Remark}\label{6.5} If $\aleph$ is any cardinal 
and $\A$ is any abelian group with torsion free
rank $\aleph$, then $\Hom(\A,\R)\cong\R^\aleph$.
\end{Remark}

\begin{Proof} In \cite{compbook}, Theorem 8.20, 
pp.387ff.~it is discussed
that $G$ contains totally disconnected compact subgroups $\Delta$ such that
the annihilator in the character group of $G$, say,
$\Delta^\perp\subseteq \A$ is free, and $\A/\Delta^\perp$ is a torsion group.
This means that $G/\Delta$ is a torus. We note that the inclusion
$\Delta^\perp\to\A$ induces an isomorphism 
$\K\otimes_\Z \Delta^\perp\to \K\otimes_\Z \A$
and the (torsion free) rank of $\A$ is $\rank \Delta^\perp$. 
If $\Delta^\perp\cong\Z^{(X)}$
for a set $X$ of cardinality $\rank \Delta^\perp$, 
then $\Hom(\A,\K)\cong \K^X$.
\end{Proof}

\subsection{The exponential function of $\C[G]=\C^\A$}

\noindent
We recall from Theorem 2.2 
that every weakly complete associative
unital algebra $W$  has an exponential function, which
is immediate in the case of $W=\C^\A$ as it is calculated componentwise.
If the weakly complete algebra $W$ is even a Hopf algebra, such as $\C^\A$,
then the group $\Gamma(W)$ of grouplike elements
is a pro-Lie group and $\Pi(W)$ is the (real!) Lie algebra of
the pro-Lie group $\Gamma(W)$
(\cite{dhtwo}, Theorem 6.15). If indeed
$W=\C^\A=\C[G]$ for a compact abelian group $G$, then
the exponential function 
$\exp_{\Gamma(W)}\colon\L(\Gamma(W))\to\Gamma(W))$
of $\Gamma(W)$ is the restriction of the (componentwise!) 
exponential function
$\exp\colon \C^\A\to (\C^\A)^{-1}=((\C^\times,.))^\A$ 
to $\L(\Gamma(W))=\Hom(\A,\C)$.

\section{The Weakly Complete Enveloping Algebra\\ 
of a Weakly Compact  Lie Algebra}

 We have observed that for compact groups $G$ the weakly complete
real group algebra $\R[G]$ contains a substantial volume of 
different materials: the pro-Lie group $G$ itself,
its Lie algebra $\L(G)$, 
the exponential function
between them and, as was discussed in  detail in \cite{dhtwo}, 
a substantial portion of the Radon measure theory of $G$. 
The topological  Hopf algebra  $\K[G]$
is, in a sense, universally generated by $G$. So it seems natural
to ask the question whether $\L(G)$ generates $\K[G]$ in a universal
way---perhaps in some fashion that would resemble the universal enveloping
algebra of a Lie algebra such as it is presented in the famous
{\sc Poincar\'e-Birkhoff-Witt}-Theorem 
(see e.g.\ \cite{bour}, \S2, n$^{\rm o}$ 7,
Th\'eor\`eme 1., p.30). This is not exactly the case,
but a few aspects  should be discussed.

So we let $\K$ again denote one of the topological fields $\R$ or $\C$.
Let $\AA$ denote the category of weakly complete associative 
unital algebras over
$\K$ and  and
$\LL$ the category of weakly complete Lie algebras over $\K$.
The functor $A\mapsto \li A$ which associates with a weakly 
complete associative algebra $A$ the weakly complete Lie algebra obtained
by considering on the weakly complete vector space $A$ the Lie
algebra obtained with respect to the Lie bracket $[x,y]=xy-yx$ is
called the {\it underlying Lie algebra functor}. Since
$A$ is a strict projective limit of finite dimensional $K$-algebras
by \cite{dhtwo}, Theorem 3.2, then  $\li A$ is a strict projective limit
of finite dimensional $\K$-Lie algebras, briefly called 
{\it pro-Lie algebras}. Every pro-Lie algebra is weakly complete.
(Caution: A comment following
Theorem 3.12 of \cite{dhtwo} exhibits an example of a weakly complete
$\K$-Lie algebra which is not a pro-Lie algebra!)

\medskip

\begin{Lemma}  \label{8.1}  
The \emph{underlying Lie algebra} functor $A\mapsto\li A$
from $\WA$ to $\LL$ has a left adjoint $\UU\colon\LL\to\AA$.
\end{Lemma}

\begin{Proof} The category $\LL$ is complete. 
(Exercise. Cf.\ Theorem A3.48 of 
\cite{compbook}, p.\ 781.) The
 ``Solution Set Condition'' (of Definition A3.59 in \cite{compbook}, 
p.\ 786) holds.
(Exercise: Cf.\ the proof of \cite{dhtwo}, 
Section 5.1 ``The solution set condition''.) 
Hence $\UU$ exists by the Adjoint Functor Existence Theorem
(i.e., Theorem A3.60 of \cite{compbook}, p.\ 786).
\end{Proof}

\medskip

In other words, 
for each weakly complete Lie algebra $L$ there is a natural morphism
$\lambda_L\colon L\to \UU(L)$ such that for each continuous Lie algebra
 morphism $f\colon L\to \li A$ for a weakly complete associative unital 
algebra $A$ there is a unique $\WA$-morphism $f'\colon \UU(L)\to A$ 
such that
$f=\li{f'}\circ\lambda_L$. 

 $$
\begin{matrix}& \LL&&\hbox to 7mm{} &\AA\cr 
\noalign{\vskip3pt}
\noalign{\hrule}\cr
\noalign{\vskip3pt}%
   L&\mapright{\lambda_L}&\li{\UU(L)}&\hbox to 7mm{} &\UU(L)\\
\lmapdown{\forall f}&&\mapdown{\li{f'}}&\hbox to 7mm{}&
         \mapdown{\exists! f'}\\
 \li A&\lmapright{\id}&\li A&\hbox to 7mm{}&A.
\end{matrix}
 $$
  
\medskip

\noindent If necessary we shall write $\UU_\K$ instead of $\UU$ 
whenever the ground field should be emphasized.
We shall call $\UU_\K(L)$ 
{\it the weakly complete enveloping algebra} of $L$ 
(over $\K$).

\medskip

\begin{Example} \label{8.2} Let $L=\K$, the smallest possible nonzero
Lie algebra over $\K$.
Then $\UU(L)=\K\<X\>$ (see \cite{dhtwo}, Definition
following Corollary 3.3), and define 
$\lambda_L\colon L\to \li{\UU(L)}$ by
$\lambda_L(t)= t{\cdot}X$. Indeed the universal property is satisfied by
\cite{dhtwo}, Corollary 3.4. Namely, let $f\colon \K\to \li A$ 
be a morphism of weakly complete Lie algebras. 
Then there is a unique morphism $f'\colon \UU(L)\to A$
such that $f'(X)=f(1)$ by \cite{dhtwo}, Corollary 3.4. 
Then $f'(t{\cdot}X)=t{\cdot}F'(X) =t{\cdot}f(1)=f(t)$. 

\smallskip
\nin Thus by Lemma 3.5 of \cite{dhtwo} and the subsequent remarks we have:

\medskip
\hbox to \hsize{\hskip1.5cm \vbox{\hsize= .8\hsize%
\nin\it The weakly complete enveloping algebra $\UU_\C(\C)$ over $\C$ of the
smallest nonzero complex Lie algebra is isomorphic to the weakly complete 
commutative
algebra $\C[[X]]^\C$ with the complex power series algebra 
$\C[[X]]\cong\C^{\N_0}$, $\N_0=\{0,1,2,\dots\}$.}}
\end{Example}

\noindent 
The size of the weakly complete enveloping algebras 
therefore is considerable.

\begin{Proposition} \label{8.3}  
The universal enveloping functor $\UU$ is multiplicative,
that is, there is a natural isomorphism 
$\alpha_{L_1L_2}\colon \UU(L_1\times L_2)\to \UU(L_1)\otimes_\W\UU(L_2)$.
\end{Proposition}

\begin{Proof} We have a natural bilinear inclusion map of weakly complete
vector spaces 
$j\colon \UU(L_1)\times \UU(L_2)\to \UU(L_1)\otimes_\W \UU(L_2)$
yielding  
$$L_1{\times}L_2\mapright{\lambda_1{\times}\lambda_2}\li{\UU(L_1)}{\times}\li{\UU(L_2)}
\mapright{j}\li{\UU(L_1)}{\otimes_\W}\li{\UU(L_2)}$$
and
$$\li{\UU(L_1)}\otimes_\W\li{\UU(L_2)}{=}\li{(\UU(L_1)\otimes_\W\UU(L_2))},$$
the composition $\alpha_0$ of which 
is a morphism of weakly complete Lie algebras.
Hence the universal property  yields a morphism 
of weakly complete associative algebras
$$\alpha\colon\UU(L_1\times L_2)\to\UU(L_1)\otimes_\W\UU(L_2)\leqno(1)$$ 
such that 
$\alpha_0=\li{\alpha}\circ\lambda_{L_1}\otimes \lambda_{L_2}$.

The functorial property of $\UU$ allows us to argue 
that each of $\UU(L_m)$, $m=1,2$
is a retract of $\UU(L_1\times L_2)$  so that we may assume 
$\UU(L_m)\subseteq\UU(L_1\times L_2)$, $m=1,2$.
Now the multiplication in $\UU(L_1\times L_2)$ 
gives rise to a continuous bilinear
map $\UU(L_1)\times\UU(L_2)\to\UU(L_1\times L_2)$, 
and then the universal property of 
the tensor product of weakly complete vector spaces yields the
morphism
$$\beta\colon \UU(L_1)\otimes_\W\UU(L_2)\to \UU(L_1\times L_2).\leqno(2)$$
Similarly to the proof of \cite{dhtwo}, Theorem 5.5 
(preceding the statement of the theorem)
we argue that $\alpha$ and $\beta$ are inverses of each other, and so 
$\alpha$ of (1) is the desired isomorphism $\alpha_{L_1L_2}$.  
\end{Proof}

\begin{Lemma} \label{li-morph} For any weakly complete unital algebra
$A$, the vector space morphism $\Delta_A\colon A\to A\otimes_\W A$,
$\Delta_A(a)=a\otimes 1+1\otimes a$
is a morphism of weakly complete 
Lie algebras  $\li A\to \li{(A\otimes_WA)}$.
\end{Lemma}
\begin{Proof} Since the functions $a\mapsto a\otimes 1,\quad 1\otimes a$
are morphisms of topological vector spaces, so is $\Delta_A$
For $y_1, y_2\in A$, write 
$z_j\defi \Delta_A(y_j)=y_j\otimes1+1\otimes y_j$
Then
 just as in the classical case, we
calculate $[\Delta_A(y_1),\Delta_A(y_2)]=$
$$[z_1,z_2]=z_1z_2-z_2z_1
=(y_1\otimes1 + 1\otimes y_1)(y_2\otimes1+1\otimes y_2)
-(y_2\otimes1 + 1\otimes y_2)(y_1\otimes1+1\otimes y_1)$$
$$=(y_1y_2\otimes1 +y_1\otimes y_2+y_2\otimes y_1+ 1\otimes y_1y_2)
-(y_2y_1\otimes1 +y_2\otimes y_1+y_1\otimes y_2+ 1\otimes y_2y_1)$$
$$=[y_1,y_2]\otimes1+1\otimes [y_1,y_2]= \Delta_A[y_1,y_2].$$
Thus $\Delta_A$ is a morphism of 
Lie algebras as asserted.
\end{Proof}

\nin 
Now we consider a weakly complete Lie algebra $L$ and
recall that $\lambda_L\colon L\to\li{U(L)}$ is a morphism of 
weakly complete Lie algebras. Thus by Lemma \ref{li-morph},

\centerline{%
 $p_L=\Delta_{U(L)}\circ \lambda_L\colon L\to \li{(\UU(L)\otimes \UU(L))}$}

\nin is a morphism of
weakly complete  Lie algebras.
Now 
by the universal property of $\UU$, $p_L$ induces
a unique natural morphism of weakly complete associative unital
algebras
$\gamma_L\colon \UU(L)\to \UU(L)\otimes_\W \UU(L)$ 
such that
$p_L=\li{(\gamma_L)}\circ \lambda_L$.
This yields the following insight, where 
$k_L\colon L\to\{0\}$ denotes the constant morphism.    

\begin{Corollary} \label{hopf}
Each weakly complete enveloping algebra $\UU(L)$ is
a weakly complete Hopf algebra with the comultiplication 
$\gamma_L$ and the coidentity $\UU(k_L)\colon\UU(L)\to\K$.
\end{Corollary}

\begin{Proof} 
 Observe that 
 $\gamma_L$ is a morphism of weakly complete unital algebras 
satisfying $\gamma_L(y)=y\otimes1+1\otimes y$ 
for $y=\lambda(x)$, $x\in L$.   
The associativity of this comultiplication is readily checked as
in the case of abstract enveloping algebras. 
The constant morphism of weakly complete 
Lie algebras $L\to\{0\}$ yields a morphism 
of weakly complete unital algebras
$\UU(L)\to \UU(\{0\})=\K$ which is the coidentity of the Hopf algebra. 
\end{Proof} 

\medskip 

\nin
Our results from \cite{dhtwo} regarding 
weakly complete associative unital algebras and
Hopf algebras over $\K$ apply to the present situation.

\medskip

\begin{Theorem} \label{8.5} {\rm(The Weakly Complete Enveloping Algebra)}
 Let $L$ be a weakly complete Lie algebra. 
Then the following statements hold:
\begin{enumerate}[\rm(i)]
\item  $\UU(L)$ is a strict projective limit 
of finite-dimensional associative unital algebras and
the group of units $\UU(L)^{-1}$ is dense in $\UU(L)$. 
It is an almost connected  pro-Lie group
(which is connected in the case of $\K=\C$).
The algebra $\UU(L)$ has an exponential function 
$\exp\colon{\li{\UU(L)}}\to \UU(L)^{-1}$,
 
\item  The pro-Lie algebra $\Pi(\UU(L))$ of primitive elements 
of $\UU(L)$ contains $\lambda_L(L)$,

\item The subalgebra generated by $\lambda_L(L)$ in $\UU(L)$ is
dense in $\UU(L)$.

\item  The pro-Lie algebra $\Pi(\UU(L))$ 
is the Lie algebra of the pro-Lie group
$\Gamma(\UU(L))$  of grouplike elements of $\UU(L)$.
We use the abbreviation  $G\defi\Gamma(\UU(L))$ and
note that  
the exponential function
$\exp_G\colon \L(G)\to G$ is 
the restriction and corestriction of the exponential
function $\exp$ of $\UU(L)$ to $\Pi(\UU(L))$,
respectively, $G$.
The image $\exp(\L(G))$ generates algebraically and topologically 
the identity component $G_0$ of $G$.
\end{enumerate}
\end{Theorem}

\begin{Proof} (i) See \cite{dhtwo}, Theorems 3.2, 3.11, 3.12, 4.1.

(ii) The very definition of the comultiplication $\gamma_L$ for 
Corollary \ref{hopf} shows that for
any $y\in \lambda_L(L)$, the image under 
the comultiplication $\gamma_L$ is
$y\otimes1+1\otimes y$, which means that $y$ is primitive.

(iii) An argument analogous to that in the proof of Proposition 5.3
of \cite{dhtwo} showing that, for the case of any topological group 
$T$, the subset  $\eta_G(T)$ of the weakly complete group
algebra $\K[T]$ spans a dense subalgebra, 
 shows here that the closed subalgebra
$S$ generated in $\UU(L)$  by $\lambda_L(L)$ has 
the universal property of $\UU(L)$
and therefore agrees with $\UU(L)$. 

(iv) See Theorem \ref{exp} and \cite{dhtwo}, Theorem 6.15.
\end{Proof}

\begin{Remark}
We note right away that for any weakly complete Lie algebra 
$L$ which has at least
one nonzero finite dimensional $\K$-linear representation, 
the morphism 
$\lambda_L\colon L\to \li{\UU(L)}$ is nonzero. By Ado's Theorem, 
this applies, in particular,
to any Lie algebra which has a nontrivial finite 
dimensional quotient and therefore
is true for the Lie algebra $\L(P)$ of any pro-Lie group $P$.
\end{Remark}

\begin{Corollary}  \label{8.6} 
{\rm (i)}  The weakly complete enveloping algebra $\UU(L)$ of a
weakly complete 
Lie algebra $L$ with a nontrivial
finite dimensional quotient has nontrivial grouplike elements. 
 
\smallskip
{\rm(ii)} If $L$ is a pro-Lie algebra, then 
$\lambda_L\colon L\to \li{\UU(L)}$
maps $L$ isomorphically onto a closed Lie subalgebra of 
the pro-Lie algebra  $\Pi(\UU(L))$ of primitive elements.
\end{Corollary} 
\begin{Proof}
 (i) If $\Pi(\UU(L))$ is nonzero, then 
$\UU(L)$ has nontrivial grouplike elements
by Theorem \ref{8.5} (iii), and by part (ii) 
of \ref{8.5},
 this is the case if $\lambda_L$
is nonzero which is the case for all $L$ satisfying the hypothesis of the 
Corollary by the remark preceding it. 

(ii) Since each finite dimensional quotient of $L$ has 
a faithful representation
by the Theorem of Ado, and since the finite dimensional 
quotients separate the
points of $L$, the morphism 
 $\lambda_L$ is injective. 
However, injective morphisms of
weakly complete vector spaces are open onto their images. 
\end{Proof}

It follows that {\it for pro-Lie algebras $L$ we may assume 
that $L$ is in fact
a  closed Lie subalgebra of primitive elements of $\UU(L)$
which  generates 
$\UU(L)$ algebraically and topologically as a weakly complete  algebra.}

It remains an open question under which circumstances
we then have in fact  $L = \P(\UU(L))$. 
In the classical setting of the discrete enveloping Hopf algebra
in characteristic 0 this is the case: see e.g.\ \cite{serre-1}, Theorem 5.4.

\bigskip

One application of the functor $\UU$ is of present interest to us. 
Recall that for a
compact group we naturally identify $G$ with the group of 
grouplike elements of $\R[G]$ (cf.\ \cite{dhtwo}, Theorems 8.7,
8.9 and 8.12),
and that $\L(G)$ may be identified with 
the pro-Lie algebra 
$\Pi(\R[G])$ of
primitive elements. (Cf.\  also Theorem \ref{exp} above.)
We  may also assume that $\L(G)$ is contained the set 
$\Pi(\UU(\L(G)))$ of primitive
elements of $\UU_\R(\L(G))$.

\begin{Theorem} \label{8.7} {\rm(i)}  Let $G$ be a 
compact group. Then there
is a natural morphism of weakly complete algebras 
$\omega_G\colon\UU_\R(\L(G))\to \R[G]$ fixing the elements
of $\L(G)$ elementwise.

{\rm(ii)}  The image of $\omega_G$ is the closed subalgebra 
$\R[G_0]$ of $\R[G]$.

{\rm(iii)}  The pro-Lie group $\Gamma(\UU_\R(\L(G)))$ is mapped onto 
$G_0=\Gamma(\R[G_0])\subseteq \R[G]$.  The connected pro-Lie group 
$\Gamma(\UU_R(\L(G)))_0$ maps surjectively onto  $G_0$ and $\P(\UU_R(\L(G)))$
onto $\P(\R[G])$.
\end{Theorem}

\begin{Proof} (i) follows at once from the universal property of $\UU$. 

(ii) As a morphism
of weakly complete Hopf algebras, $\omega_G$ has a closed image which is
generated as a weakly complete subalgebra by $\L(G)$ which is
$\R[G_0]$ by Corollary \ref{7.3} (ii). 

(iii) The morphism $\omega_G$ of weakly complete Hopf algebras maps
grouplike elements to grouplike elements, 
whence we have the commutative diagram
$$\begin{matrix}
\L(G)\subseteq\Pi(\UU_\R(\L(G)))&\mapright{\Pi(\omega_G)}&
       \Pi(\R(G)){=}\L(G)\\
\lmapdown{\exp_{\Gamma(\UU_\R(\L(G)))}} &&\mapdown{\exp_G}\\
\Gamma(\UU_\R(\L(G)))&\lmapright{\Gamma(\omega_G)}&
                         \Gamma(\R[G])=G.
\end{matrix}$$
Since $\P(\omega)$ is a retraction and the image of $\exp _G$
topologically generates $G_0$,  the image
of $\G(\omega_G)\circ \exp_{\UU_\R(\L(G))}$ topologically generates
$G_0$. Since the image of the exponential function of the pro-Lie group 
$\G(\UU_\R(\L(G)))$ generates topologically its identity
component,  $\G(\omega_G)$ maps this identity component onto $G_0$.

Since $\L(G)\subseteq \P(\UU_\R(\L(G))$, and since also any morphism
of Hopf algebras maps a primitive element onto a primitive
element we know $\omega_G(\P(\UU_\R(\L(G))))=\P(\R[G])$.\end{Proof}

\nin
It remains an open question whether $\G(\UU_\R(\L(G)))$ is in fact
connected. 

\medskip

\noindent An overview of the situation may be helpful:

\vskip-25pt
$$
\begin{matrix}%
   &&&& \R[G]\\ 
   &&&& \Big|\\
\UU_\R(\L(G))&\mapright{\omega_G,onto}&\R[G_0]&&|\\
 \Big| &&\Big|&&\Big|\\
\Gamma(\UU_\R(G))&\mapright{}&\big|&\mapright{}&\Gamma(\R[G])=G\\
\Big|&&\Big|&&\Big|\\
\Gamma(\UU_\R(\L(G)))_0&\mapright{onto}&\Gamma(\R[G_0]){=}G_0&
                              = &G_0\\
\mapup{\exp_{\UU_\R(\L(G))}}&&\mapup{\exp_{\R[G]}}&=&\mapup{\exp_G}\\
\Pi(\UU_\R(\L(G)))&=&\Pi(\R[G_0]){=}\Pi(\R[G])&=&\L(G).\\
\end{matrix}$$

\bsk

\nin
A {noteworthy} consequence of the preceding results is the insight
that {\it for any nonzero weakly complete real Lie algebra of the kind 
$L=\R^X\times \prod_{j\in J}L_j$ for any set $X$ and any family of
compact finite dimensional simple Lie algebras $L_j$, 
the weakly complete enveloping algebra $\UU(L)$ has grouplike
elements.} In the discrete situation, the enveloping algebra
 $U(L)$ of a Lie algebra $L$ for  characteristic zero has no grouplike
elements.

\bsk
    
\nin{\bf Acknowledgments.} We thank the referee for carefully reading
 our manuscript and so eliminating a good number of typos and flaws.
--
 An essential part of this
 text was written while the authors 
were partners in the program {\sc Research 
in Pairs} at the Mathematisches Forschungsinstitut Oberwolfach MFO 
in the {Black Forest} from February 4 through 23, 2019. The authors
are grateful for  the environment 
and infrastructure of MFO which made this research possible. 
--
The authors
also acknowledge the joint work with {\sc Rafael Dahmen} of the
Karlsruhe Institut f\"ur Technologie;  a good deal of that work 
did take place 
at the Technische Universit\"at Darmstadt in 2018).
--
{The second author thanks Siegfried Echterhoff for helpful insights.}

\end{document}